\documentclass[12pt, reqno]{amsart}

\usepackage{times}
\usepackage{amsfonts,amsmath}
\usepackage{amssymb}
\usepackage[colorlinks, linkcolor=blue, anchorcolor=blue, citecolor=blue]{hyperref}

\usepackage{geometry}
\newtheorem*{thma}{Theorem A}
\newtheorem*{thmb}{Theorem B}

\newtheorem*{declaration}{Declaration}

\newtheorem{theorem}{Theorem}[section]
\newtheorem{corollary}[theorem]{Corollary}
\newtheorem{lemma}[theorem]{Lemma}
\newtheorem{proposition}[theorem]{Proposition}

\theoremstyle{definition}

\theoremstyle{remark}
\newtheorem{remark}[theorem]{Remark}

\numberwithin{equation}{section}

\allowdisplaybreaks

\begin{document}

\title[Generalized Volterra-type integral operators between Bloch-type spaces]
{Generalized Volterra-type integral operators between Bloch-type spaces}

\author[C. Tong]{Cezhong Tong}
\address[C. Tong]{Department of Mathematics, Hebei University of Technology, Tianjini 300401, China.}
\email{ctong@hebut.edu.cn, cezhongtong@hotmail.com}

\author[X. He]{Xin He}
\address[X. He]{Department of Mathematics, Hebei University of Technology, Tianjini 300401, China.}
\email{hexinsx2000@163.com}

\author[Z. Yang]{Zicong Yang}
\address[Z. Yang]{Department of Mathematics, Hebei University of Technology, Tianjin 300401, China.}
\email{zicongyang@126.com}

\thanks{The authors were supported in part by the National Natural Science Foundation of China (Grant No. 12171136), the Natural Science Foundation of Hebei Province (Grant No. A2020202005, A2023202031, A2023202037), and the Natural Science Foundation of Tianjin City (Grant No. 20JCYBJC00750).}

\begin{abstract}
The Volterra-type integral operator plays an essential role in modern complex analysis and operator theory. Recently, Chalmoukis \cite{Cn} introduced a generalized integral operator, say $I_{g,a}$, defined by
$$I_{g,a}f=I^n(a_0f^{(n-1)}g'+a_1f^{(n-2)}g''+\cdots+a_{n-1}fg^{(n)}),$$
where $g\in H(\mathbb{D})$ and $a=(a_0,a_1,\cdots,a_{n-1})\in \mathbb{C}^n$. $I^n$ is the $n$th iteration of the integral operator $I$. In this paper, we introduce a more generalized integral operators $I_{\mathbf{g}}^{(n)}$ that cover $I_{g,a}$ on the Bloch-type space $\mathcal{B}^{\alpha}$, defined by 
$$I_{\mathbf{g}}^{(n)}f=I^n(fg_0+\cdots+f^{(n-1)}g_{n-1}).$$
We show the rigidity of the operator $I_{\mathbf{g}}^{(n)}$ and further the sum $\sum_{i=1}^nI_{g_i}^{N_i,k_i}$, where $I_{g_i}^{N_i,k_i}f=I^{N_i}(f^{(k_i)}g_i)$. Specifically, the boundedness and compactness of $\sum_{i=1}^nI_{g_i}^{N_i,k_i}$ are equal to those of each $I_{g_i}^{N_i,k_i}$. Moreover, the boundedness and compactness of $I^n((fg')^{(n-1)})$ are independent of $n$ when $\alpha>1$.
\end{abstract}

\maketitle

\section{Introduction}

Let $\mathbb{D}$ be the open unit disk in the complex plane $\mathbb{C}$ and $H(\mathbb{D})$ be the space of all analytic functions on $\mathbb{D}$. For $\alpha>0$, a function $f\in H(\mathbb{D})$ is called an $\alpha$-Bloch function if 
$$b_{\alpha}(f)=\sup\{(1-|z|^2)^{\alpha}|f'(z)|:z\in\mathbb{D}\}<\infty.$$
This defines a seminorm, and the family of all $\alpha$-Bloch functions form a complex Banach space $\mathcal{B}^{\alpha}$ with the norm
$$\|f\|_{\mathcal{B}^{\alpha}}=|f(0)|+b_{\alpha}(f).$$
The little $\alpha$-Bloch space $\mathcal{B}_0^{\alpha}$ is a closed subspace of $\mathcal{B}^{\alpha}$ and consists of all $f\in H(\mathbb{D})$ such that
$$\lim_{|z|\to 1^-}(1-|z|^2)^{\alpha}|f'(z)|=0.$$
When $\alpha=1$,  we abtain the Bloch functions and corresponding Bloch space, which is denoted by $\mathcal{B}$. For the general theory of $\alpha$-Bloch functions, see \cite{LbOyc, Zkh}.

Let $g\in H(\mathbb{D})$, the volterra-type integral operators on $H(\mathbb{D})$ is defined by 
\begin{equation*}
I_gf(z)=\int_{0}^zf(\zeta)g'(\zeta)d\zeta\quad {\rm and}\quad J_gf(z)=\int_{0}^zf'(\zeta)g(\zeta)d\zeta.
\end{equation*}
The importance of the operators $I_g$ and $J_g$ comes from the fact that 
$$I_gf(z)+J_gf(z)=M_gf(z)-g(0)f(0),$$
where $M_gf=g\cdot f$ is the multiplication operator.  The operator $I_g$ is also seen as a generalization of the classical Ces\`aro averaging operator. After the original work of Pommerenke \cite{Pc}, the operator $I_g$ has been studied extensively on various analytic function spaces. See \cite{AaCja, AaSag} for Hardy spaces, \cite{AaCo, AaSag1, Xj1,YyLj} for Bergman spaces, \cite{GpGdPja} for Dirichlet spaces, \cite{FzsZzh, Hzj, Xj} for Bloch-type spaces and \cite{Co, CoPj, Us} for Fock spaces.

Let $n$ be a positive integer and $a=(a_0, a_1, \cdots, a_{n-1})\in \mathbb{C}^n$. Recently, Chalmoukis \cite{Cn} introduced a  class of generalized integral operator $I_{g,a}$ on Hardy spaces, which is defined as follows:
$$I_{g,a}f=I^n\left(a_0fg^{(n)}+a_1f'g^{(n-1)}+\cdots+a_{n-1}f^{(n-1)}g'\right),$$
where $If(z)=\int_{0}^zf(\zeta)d\zeta$ is the integration operator and $I^n$ is the $n$th iteration of $I$. The operator $I_{g,a}$ is clearly a generalization of $I_{g}$ since 
$$I_{g,a}f=I^n((fg')^{n-1})$$
when $a_k=\binom{n-1}{k}$ for $k=0,1,\cdots,n-1$. Chalmoukis completely characterized the boundedness and compactness of $I_{g,a}:H^p\to H^q$ when $0<p\leq q<\infty$, here $H^p$ is the classical Hardy space over $\mathbb{D}$. 

Instead of a single function, we consider an analytic vector-valued function $\mathbf{g}$ on $\mathbb{D}$, namely, $\mathbf{g}=(g_0, g_1, \cdots, g_{n-1})$ where $g_k\in H(\mathbb{D})$ for $k=0,1,\cdots,n-1$. We define a more generalized integral operator $I_{\mathbf{g}}^{(n)}$ on $H(\mathbb{D})$ as follows:
$$I_{\mathbf{g}}^{(n)}f=I^n\left(fg_0+f'g_1+\cdots+f^{(n-1)}g_{n-1}\right).$$
It is clear that $I_{\mathbf{g}}^{(n)}=\sum_{k=0}^{n-1}I_{g_k}^{n,k}$, where for $g\in H(\mathbb{D})$, $I_{g}^{n,k}f=I^n(f^{(k)}g)$. And $I_{\mathbf{g}}^{(n)}$ coincides with Chalmoukis' integral operator $I_{g,a}$ when $g_k=a_kg^{(n-k)}$ for some $g\in H(\mathbb{D})$ and $k=0,1,\cdots, n-1$.

In this paper, we will show the rigidity of the operator $I_{\mathbf{g}}^{(n)}$ and then $I_{g,a}$ from $\mathcal{B}^{\alpha}$ to $\mathcal{B}^{\beta}$. Our main results are stated as follows.      

\begin{thma}\label{theoremA}
Suppose $\alpha, \beta>0$,  and $\mathbf{g}=(g_0, g_1, \cdots, g_{n-1})$ with $g_k\in H(\mathbb{D})$ for $0\leq k\leq n-1$. Then $I_{\mathbf{g}}^{(n)}:\mathcal{B}^{\alpha}\to \mathcal{B}^{\beta}$ is bounded (or compact, resp.) if and only each $I_{g_k}^{n,k}:\mathcal{B}^{\alpha}\to \mathcal{B}^{\beta}$, $(0\leq k\leq n-1)$, is bounded (or compact, resp.)
\end{thma}

\begin{thmb}\label{theoremB}
Let $\alpha,\beta>0$ and $g_i\in H(\mathbb{D})$ for $i=1,2,\cdots,n$. Suppose $N_i,k_i\in\mathbb{N}$ with $0<t_i=N_i-k_i<N_i$ and $t_i$ are all different. Then $\sum_{i=1}^n I_{g_i}^{N_i,k_i}:\mathcal{B}^{\alpha}\to \mathcal{B}^{\beta}$ is bounded (or compact, resp.) if and only if each $I_{g_i}^{N_i,k_i}:\mathcal{B}^{\alpha}\to\mathcal{B}^{\beta}$, $1\leq i\leq n$, is bounded (or compact, resp.).
\end{thmb}

The paper is organized as follows. In Section 2, we collect some well-known facts and prove a key lemma. We give the proof of Theorem A and Theorem B in Section 3. Furthermore, we can give an intact sufficient and necessary conditions for the boundedness and compactness of $I_{g}^{n,k}$ in the case ``$\alpha>1$" or ``$\alpha>0$ and $k\neq 0$". The case ``$k=0$ and $0<\alpha\leq 1$" is studied in Section 4. In the end of the paper, we additionally obtain conditions for certain linear differential equations to have solutions in Bloch-type spaces.

Throughout the paper, we write $A\lesssim B$ (or $B\gtrsim A$) if there exists an absolute constant $C>0$ such that $A\leq CB$. As usual, $A\simeq B$ means $A\lesssim B$ and $B\lesssim A$. We will be more specific if the dependence of such constants on certain parameters becomes critical.
 
\section{Preliminries}

In this section, we present some preliminary facts and auxiliary lemmas which will be needed in our paper. We begin with the following well-known characterization for Bloch-type space via higher order derivatives. This result has been obtained in \cite[Proposition 8]{Zkh}. See also \cite[Proposition 2.1]{LbOyc}

\begin{lemma}\label{lemma2.1}
Suppose $\alpha>0$ and $n$ is a positive integer. An analytic function $f$ belongs to $\mathcal{B}^{\alpha}$ if and only if 
\begin{equation*} 
\sup_{z\in\mathbb{D}}(1-|z|^2)^{\alpha+n}|f^{(n+1)}(z)|<\infty.
\end{equation*}
And $f\in \mathcal{B}_{0}^{\alpha}$ if and only if 
\begin{equation*}
\lim_{|z|\to 1^-}(1-|z|^2)^{\alpha+n}|f^{(n+1)}(z)|<\infty.
\end{equation*}
Moreover,
\begin{equation*}
\|f\|_{\mathcal{B}^{\alpha}}\simeq \sum_{j=0}^{n}|f^{(j)}(0)|+\sup_{z\in\mathbb{D}}(1-|z|^2)^{\alpha+n}|f^{(n+1)}(z)|.
\end{equation*}
\end{lemma}

\begin{lemma}[\cite{Ss}]\label{lemma2.2}
Suppose $\alpha>0$ and $f\in\mathcal{B}^{\alpha}$, then for all $z\in\mathbb{D}$,
\begin{equation*}
|f(z)|\leq C\left\{
\begin{aligned}
&|f(0)|+b_{\alpha}(f),  &0<\alpha<1;\\
&|f(0)|+b_{\alpha}(f)\log\frac{2}{1-|z|^2}, & \alpha=1;\\
&|f(0)|+b_{\alpha}(f){(1-|z|^2)^{\alpha-1}}, &\alpha>1
\end{aligned}
\right.
\end{equation*}
for some constant $C>0$ independent of $f$.
\end{lemma}

Lemma \ref{lemma2.2} tells us that the evaluation functions are bounded in $\mathcal{B}^{\alpha}$. Applying Montel's Theorem of normal family, we have the following lemma. The proof is just a modification of \cite[Proposition 3.11]{CcMbd}.

\begin{lemma}\label{lemma2.3}
The operator $I_{g}^{n,k}:\mathcal{B}^{\alpha}\to \mathcal{B}^{\beta}$ is compact if and only if $I_{g}^{n,k}$ is bounded and $\|I_{g}^{n,k}f_j\|_{\alpha}\to 0$ for any bounded sequence $\{f_j\}_{j=1}^{\infty}$ in $\mathcal{B}^{\alpha}$ which converges to zero uniformly on compact subsets of $\mathbb{D}$. 
\end{lemma}

For any $w\in\mathbb{D}$ and $i\in\mathbb{N}$, let 
$$
f_{w}^{[i]}(z)=\frac{z^i}{(1-\overline{w}z)^{\alpha+i}},\quad z\in\mathbb{D}.
$$ 
Then we could get the following estimates for the norm of $f_w^{[i]}$.

\begin{lemma}\label{lemma2.4}
$$\|f_w^{[i]}\|_{\alpha}\simeq \frac{1}{(1-|w|^2)^{i+1}}$$
for all $w\in\mathbb{D}$ and $i\geq 1$. And 
$$\|f_w^{[0]}\|_{\alpha}\simeq \frac{1}{1-|w|^2}$$
for all $w\in\mathbb{D}$ with $|w|>\frac{1}{2}$.
\end{lemma}

\begin{proof}
On one hand, note that $|1-\overline{w}z|\geq \min\{1-|z|, 1-|w|\}$, then we have 
\begin{equation}\label{equa2.1}
\begin{split}
\|f_w^{[i]}\|_{\mathcal{B}^{\alpha}}&=|f_w^{[i]}(0)|+\sup_{z\in\mathbb{D}}(1-|z|^2)^{\alpha}\left|(f_{w}^{[i]})'(z)\right|\\
&\lesssim 1+\sup_{z\in\mathbb{D}}(1-|z|^2)^{\alpha}\frac{1}{|1-\overline{w}z|^{\alpha+i+1}}\\
&\lesssim \frac{1}{(1-|w|^2)^{i+1}}.
\end{split}
\end{equation}
On the other hand, when $|w|>\frac{1}{2}$, 
\begin{equation*}
\|f_w^{[0]}\|_{\mathcal{B}^{\alpha}}\geq(1-|w|^2)^{\alpha}\left|(f_w^{[0]})'(w)\right|\geq \frac{\alpha}{2}\frac{1}{1-|w|^2}.
\end{equation*}
Thus, $\|f_w^{[0]}\|_{\mathcal{B}^{\alpha}}\simeq \frac{1}{1-|w|^2}$ for all $w\in \mathbb{D}$ with $|w|>\frac{1}{2}$.

When $i\geq 1$, note that $(f_w^{[i]})^{(j)}(0)=0$ for $0\leq j\leq i-1$, then by Lemma \ref{lemma2.1}, we have
\begin{equation*}
\begin{split}
\|f_{w}^{[i]}\|_{\mathcal{B}^{\alpha}}&\simeq \sup_{z\in\mathbb{D}}(1-|z|^2)^{\alpha+i-1}\left|(f_w^{[i]})^{(i)}(z)\right|\\
&\geq (1-|w|^2)^{\alpha+i-1}\left|(f_w^{[i]})^{(i)}(w)\right|\\
&\gtrsim (1-|w|^2)^{\alpha+i-1}\sum_{j=0}^{i}\binom{i}{j}\frac{|w|^{2j}}{(1-|w|^2)^{\alpha+i+j}}\\
&=\frac{1}{(1-|w|^2)^{i+1}}.
\end{split}
\end{equation*}
This, together with \eqref{equa2.1}, shows that $\|f_{w}^{[i]}\|_{\mathcal{B}^{\alpha}}\simeq \frac{1}{(1-|w|^2)^{i+1}}$ for all $w\in\mathbb{D}$.
\end{proof}

The following lemma illustrates the key idea of our main results. We believe that this result should be of some independent value and interest.

\begin{lemma}\label{lemma2.5}
Let $c_0, c_1,\cdots, c_n$ be non-zero complex scalars, then 
\begin{equation}\label{equa2.2}
\left\|\sum_{i=0}^nc_if_w^{[i]}\right\|_{\mathcal{B}^{\alpha}}\gtrsim\sum_{i=0}^n|c_i|\|f_w^{[i]}\|_{\mathcal{B}^{\alpha}}
\end{equation}
for all $w\in\mathbb{D}$ with $|w|>\frac{1}{2}$.
\end{lemma}

\begin{proof}
For $j\in\mathbb{N}$, let 
$$r_j(w)=(1-|w|^2)^{\alpha+j}\sum_{i=0}^{n}c_i\left(f_w^{[i]}\right)^{(j+1)}(w).$$ 
It follows from lemma \ref{lemma2.1} that $\left\|\sum_{i=0}^nc_if_w^{[i]}\right\|_{\mathcal{B}^{\alpha}}\gtrsim |r_j(w)|$.

The case $n=0$ in \eqref{equa2.2} is trivial. Now we prove the case $n=1$. Note that
\begin{equation*}
|r_0(w)|=\left|c_0\frac{\alpha\overline{w}}{1-|w|^2}+c_1\frac{1+\alpha|w|^2}{(1-|w|^2)^2}\right|,
\end{equation*}
and 
\begin{equation*}
|r_1(w)|=\left|c_0\frac{\alpha(\alpha+1)\overline{w}^2}{1-|w|^2}+c_1\frac{2(\alpha+1)\overline{w}+\alpha(\alpha+1)\overline{w}|w|^2}{(1-|w|^2)^2}\right|.
\end{equation*}
By triangle inequality, 
\begin{equation*}
|c_1|\frac{(\alpha+1)|w|}{(1-|w|^2)^2}\leq |r_1(w)|+(\alpha+1)|\overline{w}r_0(w)|\lesssim \left\|c_0f_w^{[0]}+c_1f_w^{[1]}\right\|_{\mathcal{B}^{\alpha}}.
\end{equation*}
When $|w|>\frac{1}{2}$, it follows from Lemma \ref{lemma2.4} that 
$$|c_1|\|f_w^{[1]}\|_{\mathcal{B}^{\alpha}}\lesssim \left\|c_0f_w^{[0]}+c_1f_w^{[1]}\right\|_{\mathcal{B}^{\alpha}}.$$
Then we also have 
$$|c_0|\|f_w^{[0]}\|_{\mathcal{B}^{\alpha}}\leq \left\|c_0f_w^{[0]}+c_1f_w^{[1]}\right\|_{\mathcal{B}^{\alpha}}+|c_1|\|f_w^{[1]}\|_{\mathcal{B}^{\alpha}}\lesssim\left\|c_0f_w^{[0]}+c_1f_w^{[1]}\right\|_{\mathcal{B}^{\alpha}}.$$

For the general case $n\geq 2$, we only need to prove 
\begin{equation}\label{equa2.3}
\left\|\sum_{i=0}^nc_if_w^{[i]}\right\|_{\mathcal{B}^{\alpha}}\gtrsim |c_n|\|f_w^{[n]}\|_{\mathcal{B}^{\alpha}}
\end{equation}
After achieving \eqref{equa2.3}, we can see
\begin{equation*}
\left\|\sum_{i=0}^{n-1}c_if_w^{[i]}\right\|_{\mathcal{B}^{\alpha}}\leq \left\|\sum_{i=0}^nc_if_w^{[i]}\right\|_{\mathcal{B}^{\alpha}}+|c_n|\|f_w^{[n]}\|_{\mathcal{B}^{\alpha}}\lesssim\left\|\sum_{i=0}^nc_if_w^{[i]}\right\|_{\mathcal{B}^{\alpha}}.
\end{equation*}
Then the inductions works to get
\begin{equation*}
|c_j|\|f_w^{[j]}\|_{\mathcal{B}^{\alpha}}\lesssim\left\|\sum_{i=0}^{n-1}c_if_w^{[i]}\right\|_{\mathcal{B}^{\alpha}}\lesssim\left\|\sum_{i=0}^nc_if_w^{[i]}\right\|_{\mathcal{B}^{\alpha}}
\end{equation*}
for $j=n-1,\cdots,1,0$. In fact, note that $(f_w^{[i]})'(z)=\frac{h_i(z)}{(1-\overline{w}z)^{\alpha+n}}$, where 
\begin{equation*}
h_i(z)=\left\{
\begin{aligned}
&\alpha\overline{w}(1-\overline{w}z)^{n-1},  &i=0;\\
&(iz^{i-1}+\alpha\overline{w}z^i)(1-\overline{w}z)^{n-i-1}, & 1\leq i\leq n-1;\\
&\frac{nz^{n-1}+\alpha\overline{w}z^n}{1-\overline{w}z}, &i=n.
\end{aligned}
\right.
\end{equation*}
Then by \cite[Lemma2.7]{AhLhjTczYzc} and triangle inequality,
\begin{equation*}
\begin{split}
&\quad r_n(w)+\sum_{i=0}^{n-1}(-1)^{i+1}\binom{n}{i+1}\Gamma_{\alpha+1}^{n-1-i,n-1}\overline{w}^{i+1}r_i(w)\\
&=c_n h_n^{(n)}(w)=c_n \frac{n!(n+\alpha)w}{(1-|w|^2)^{n+1}},
\end{split}
\end{equation*}
where $\Gamma_{\alpha+1}^{n-1-i,n-1}=(\alpha+n-i)\cdots(\alpha+n)$. When $|w|>\frac{1}{2}$, it follows from Lemma \ref{lemma2.4} that
$$|c_n|\|f_w^{[n]}\|_{\mathcal{B}^{\alpha}}\lesssim |c_n|\frac{n!(n+\alpha)|w|}{(1-|w|^2)^{n+1}}\lesssim \sum_{i=0}^n|r_i(w)|\lesssim\left\|\sum_{i=0}c_if_w^{[i]}\right\|_{\mathcal{B}^{\alpha}}.$$
Now the proof is complete.
\end{proof}

\section{Proof of Theorem A and Theorem B}
In this section, we give the proof of our main results. To this end, we need to characterize the boundedness and compactness of $T_{g}^{n,k}$ for $n\geq 1$ and $0\leq k<n$. The case ``$k=0$ and $0<\alpha\leq 1$" is different, which will be studied in the next section.

\begin{proposition}\label{proposition3.1}
Suppose $\alpha,\beta>0$, and $0<k\leq n-1$, then
\begin{itemize}
\item[(i)] $I_{g}^{n,k}:\mathcal{B}^{\alpha}\to\mathcal{B}^{\beta}$  is bounded if and only if 
\begin{equation}\label{equa3.1}
\sup_{z\in\mathbb{D}}(1-|z|^2)^{n-k+\beta-\alpha}|g(z)|<\infty.
\end{equation}
\item[(ii)] $I_{g}^{n,k}:\mathcal{B}^{\alpha}\to\mathcal{B}^{\beta}$ is compact if and only if 
\begin{equation}\label{equa3.2}
\lim_{|z|\to 1^{-}}(1-|z|^2)^{n-k+\beta-\alpha}|g(z)|=0.
\end{equation}
\end{itemize}
\end{proposition}

\begin{proof}
We begin with the proof of the ``only if part" of (i). Assume $I_{g}^{n,k}:\mathcal{B}^{\alpha}\to \mathcal{B}^{\beta}$ is bounded, then by Lemma \ref{lemma2.4},
\begin{equation}\label{equa3.3}
\|I_{g}^{n,k}f_w^{[k]}\|_{\beta}\lesssim \|f_w^{[k]}\|_{\alpha}\simeq \frac{1}{(1-|w|^2)^{k+1}}
\end{equation}
for all $w\in\mathbb{D}$. Notice that $(I_{g}^{n,k}f_w^{[k]})^{(j)}(0)=0$ for $0\leq j\leq n-1$. We apply Lemma \ref{lemma2.1} to obtain
\begin{equation*}
\begin{split}
\|I_{g}^{n,k}f_w^{[k]}\|_{\beta}&\simeq \sup_{z\in\mathbb{D}}(1-|z|^2)^{\beta+n-1}\left|(I_{g}^{n,k}f_w^{[k]})^{(n)}(z)\right|\\
&=\sup_{z\in\mathbb{D}}(1-|z|^2)^{\beta+n-1}\left|(f_w^{[k]})^{(k)}(z)\right||g(z)|\\
&\geq (1-|w|^2)^{\beta+n-1}\left|(f_w^{[k]})^{(k)}(w)\right||g(w)|\\
&\gtrsim (1-|w|^2)^{\beta-\alpha+n-2k-1}|g(w)|.
\end{split}
\end{equation*}
This, together with \eqref{equa3.3}, shows that
$$\sup_{w\in\mathbb{D}}|g(w)|(1-|w|^2)^{\beta-\alpha+n-k}<\infty.$$

Conversely, assume \eqref{equa3.1} holds. Then for any $f\in\mathcal{B}^{\alpha}$, we use Lemma \ref{lemma2.1} to get
\begin{equation*}
\begin{split}
\|I_{g}^{n,k}f\|_{\beta}&\simeq \sup_{z\in\mathbb{D}}(1-|z|^2)^{\beta+n-1}\left|(I_{g}^{n,k}f)^{(n)}(z)\right|\\
&=\sup_{z\in\mathbb{D}}(1-|z|^2)^{\beta+n-1}|f^{(k)}(z)||g(z)|\\
&\leq \sup_{z\in\mathbb{D}}(1-|z|^2)^{\beta-\alpha+n-k}|g(z)|\sup_{z\in\mathbb{D}}(1-|z|^2)^{\alpha+k-1}|f^{(k)}(z)|\\
&\lesssim \sup_{z\in\mathbb{D}}(1-|z|^2)^{\beta-\alpha+n-k}|g(z)|\|f\|_{\alpha},
\end{split}
\end{equation*}
which clarifies the boundedness of $I_{g}^{n,k}$.

The proof for the ``only if part" of (ii) is similar to (i). Indeed, for any $w\in\mathbb{D}$, let
$$q_w^{[0]}(z)=(1-|w|^2)f_w^{[0]}(z)=\frac{1-|w|^2}{(1-\overline{w}z)^{\alpha}},\quad z\in\mathbb{D}.$$
It follows from Lemma \ref{lemma2.4} that $q_w^{[0]}\in\mathcal{B}^{\alpha}$ and $\sup_{w\in\mathbb{D}}\|q_w^{[0]}\|_{\alpha}\lesssim 1$. It is also clear that $\{q_w^{[0]}\}$ converges to 0 uniformly on compact subsets of $\mathbb{D}$ as $|w|\to 1^{-}$. If $I_{g}^{n,k}:\mathcal{B}^{\alpha}\to \mathcal{B}^{\beta}$ is compact, then by Lemma \ref{lemma2.3},
\begin{equation}\label{equa3.4}
\lim_{|w|\to 1^{-}}\|I_{g}^{n,k}q_w^{[0]}\|_{\mathcal{B}^{\beta}}=0.
\end{equation}
Another application of Lemma \ref{lemma2.1} shows that
\begin{equation*}
\begin{split}
\|I_{g}^{n,k}q_w^{[0]}\|_{\mathcal{B}^{\beta}}&\simeq \sup_{z\in\mathbb{D}}(1-|z|^2)^{\beta+n-1}\left|(I_{g}^{n,k}q_w^{[0]})^{(n)}(z)\right|\\
&=\sup_{z\in\mathbb{D}}(1-|z|^2)^{\beta+n-1}|(q_w^{[0]})^{(k)}(z)||g(z)|\\
&\geq (1-|w|^2)^{\beta+n-1}|(q_w^{[0]})^{(k)}(w)||g(w)|\\
&\gtrsim (1-|w|^2)^{\beta-\alpha+n-k}|w|^k|g(w)|.
\end{split}
\end{equation*}
Thus, $\lim_{|w|\to 1^-}(1-|w|^2)^{\beta-\alpha+n-k}|g(w)|=0$.

It remains to prove the ``if part" of (ii). Let $\{f_j\}_{j=1}^{\infty}$ be an arbitrary sequence in $\mathcal{B}^{\alpha}$
such that $\sup_{j}\|f_j\|_{\mathcal{B}^{\alpha}}:=L<\infty$ and $f_j\to 0$ uniformly on compact subsets of $\mathbb{D}$. For any $\varepsilon>0$, choose $\delta\in (0,1)$ such that
\begin{equation*}
|g(z)|(1-|z|^2)^{\beta-\alpha+n-k}<\frac{\varepsilon}{L}
\end{equation*}
whenever $\delta<|z|<1$. Then we use Lemma \ref{lemma2.1} to get
\begin{equation}\label{equa3.5}
\begin{split}
\|I_{g}^{n,k}f_j\|_{\mathcal{B}^{\beta}}&\simeq \sup_{z\in\mathbb{D}}(1-|z|^2)^{\beta+n-1}\left|(I_{g}^{n,k}f_j)^{(n)}(z)\right|\\
&=\sup_{z\in\mathbb{D}}(1-|z|^2)^{\beta+n-1}|f_j^{(k)}(z)||g(z)|\\
&\lesssim \sup_{|z|\leq \delta}(1-|z|^2)^{\beta+n-1}|f_j^{(k)}(z)||g(z)|\\
&\qquad\qquad+\sup_{\delta<|z|<1}(1-|z|^2)^{\beta-\alpha+n-k}|g(z)|\|f_j\|_{\mathcal{B}^{\alpha}}\\
&\lesssim \sup_{|z|\leq \delta}(1-|z|^2)^{\beta-\alpha+n-k}|g(z)|\sup_{|z|\leq \delta}|f_j^{(k)}(z)|+\varepsilon.
\end{split}
\end{equation}
The uniform convergence of $\{f_j\}$ on compact subsets of $\mathbb{D}$ along with the Cauchy's estimate implies that $\{f_j^{(k)}\}$ also converges to 0 uniformly on compact subsets of $\mathbb{D}$. Hence
$$\lim_{j\to\infty}\sup_{|z|\leq \delta}|f_j^{(k)}(z)|=0.$$
Letting $j\to\infty$ in \eqref{equa3.5}, we get
$$\limsup_{j\to\infty}\|I_g^{n,k}f_j\|_{\alpha}\lesssim \varepsilon$$
for any $\varepsilon >0$. Consequently, $I_{g}^{n,k}:\mathcal{B}^{\alpha}\to\mathcal{B}^{\beta}$ is compact by Lemma \ref{lemma2.3}.
\end{proof}

\begin{remark}\label{remark3.2}
When $\alpha>1$, the argument in Proposition 3.1 is also applicable for the boundedness and compactness of $I_g^{n,0}$. As a consequence, $I_{g}^{n,0}:\mathcal{B}^{\alpha}\to \mathcal{B}^{\beta}$ is bounded if and only if 
$$\sup_{z\in\mathbb{D}}(1-|z|^2)^{n+\beta-\alpha}|g(z)|<\infty.$$
And $I_{g}^{n,0}:\mathcal{B}^{\alpha}\to \mathcal{B}^{\beta}$ is compact if and only if
$$\lim_{|z|\to 1^-}(1-|z|^2)^{n+\beta-\alpha}|g(z)|=0.$$
However the case ``$k=0$ and $0<\alpha\leq 1$" is different, which will be studied in the next section.
\end{remark}

We are now ready to prove our main results.

\begin{proof}[\bf Proof of Theorem A] 
The sufficiency is trivial since $I_{\mathbf{g}}^{(n)}=\sum_{k=0}^{n-1}I_{g_k}^{n,k}$. 

Now we prove the necessity. Without loss of generality, assume $n\geq 2$. For any $w\in\mathbb{D}$ and $i=1,2,\cdots,n$, let 
$$q_w^{[i]}(z)=(1-|w|^2)^{i+1}f_w^{[i]}(z)=\frac{(1-|w|^2)^{i+1}z^i}{(1-\overline{w}z)^{\alpha+1}}, \quad z\in\mathbb{D}.$$
By Lemma \ref{lemma2.4}, $q_w^{[i]}\in\mathcal{B}^{\alpha}$ and $\|q_w^{[i]}\|_{\mathcal{B}^{\alpha}}\simeq 1$. If $I_{\mathbf{g}}^{(n)}:\mathcal{B}^{\alpha}\to\mathcal{B}^{\beta}$ is bounded, then
\begin{equation}\label{equa3.6}
\sup_{w\in\mathbb{D}}\|I_{\mathbf{g}}^{(n)}q_w^{[i]}\|_{\mathcal{B}^{\beta}}<\infty,\quad i=1,2,\cdots,n.
\end{equation}
By Lemma \ref{lemma2.1}, we have
\begin{equation}\label{equa3.7}
\begin{split}
\|I_{\mathbf{g}}^nq_w^{[i]}\|_{\mathcal{B}^{\beta}}&\simeq \sup_{z\in\mathbb{D}}(1-|z|^2)^{\beta+n-1}\left|(I_{\mathbf{g}}^nq_w^{[i]})^{(n)}(z)\right|+\sum_{j=0}^{n-1}\left|(I_{\mathbf{g}}^nq_w^{[i]})^{(j)}(0)\right|\\
&=\sup_{z\in\mathbb{D}}(1-|z|^2)^{\beta+n-1}\left|\sum_{k=0}^{n-1}(q_w^{[i]})^{(k)}(z)g_k(z)\right|\\
&\geq (1-|w|^2)^{\beta+n+i}\left|\sum_{k=0}^{n-1}(f_w^{[i]})^{(k)}(w)g_k(w)\right|.
\end{split}
\end{equation}
For $i\geq 1$, 
$$f_w^{[i]}(w)=\frac{w^i}{(1-|w|^2)^{\alpha+i}}=\frac{1}{\alpha\cdots(\alpha+i-1)}\overline{(f_w^{[0]})^{(i)}(w)}.$$
Let $K(w,z)=f_w^{[0]}(z)$, when $k\geq 1$, it is easy to verify that
\begin{equation*}
(f_w^{[i]})^{(k)}(z)=\frac{1}{\alpha\cdots(\alpha+i-1)}\frac{\partial ^{k+i}K}{\partial z^k\partial \overline{w}^i}(w,z)
\end{equation*}
and
\begin{equation*}
\overline{(f_w^{[k]})^{(i)}(z)}=\frac{1}{\alpha\cdots(\alpha+k-1)}\frac{\partial ^{k+i}\overline{K}}{\partial w^k\partial \overline{z}^i}(w,z)
\end{equation*}
Therefore, $(f_w^{[i]})^{(k)}(w)=\frac{\alpha\cdots(\alpha+k-1)}{\alpha\cdots(\alpha+i-1)}\overline{(f_w^{[k]})^{(i)}(w)}$ for $i,k\geq 1$. Set 
$$\gamma_i(w)=(1-|w|^2)^{\beta+n+i}\sum_{k=0}^{n-1}(f_w^{[k]})^{(i)}(w)c_k(w),\quad i=1,2,\cdots,n,$$
where $c_0(w)=\overline{g_0(w)}$, and $c_k(w)=\alpha\cdots(\alpha+k-1)\overline{g_k(w)}$ for $1\leq k\leq n-1$. Modifying the proof of Lemma \ref{lemma2.5}, we obtain that
\begin{equation*}
\begin{split}
&\quad\gamma_n(w)+\sum_{j=1}^{n-1}(-1)^{j}\binom{n-1}{j}\Gamma_{\alpha+1}^{n-1-j,n-2}\overline{w}^j\gamma_j(w)\\
&=\frac{(n-1)!(n+\alpha-1)w}{(1-|w|^2)^n}(1-|w|^2)^{\beta-\alpha+n+1}c_{n-1}(w)\\
&=\alpha\cdots(\alpha+n-2)(\alpha+n-1)(n-1)!w(1-|w|^2)^{\beta-\alpha+1}\overline{g_{n-1}(w)}.
\end{split}
\end{equation*}
This, together with \eqref{equa3.6} and \eqref{equa3.7}, shows that
$$\sup_{w\in\mathbb{D}}(1-|w|^2)^{\beta-\alpha+1}|w||g_{n-1}(w)|\lesssim \sum_{j=1}^n|\gamma_j(w)|\lesssim \sum_{j=1}^n\|I_{\mathbf{g}}^{(n)}q_w^{[j]}\|_{\mathcal{B}^{\beta}}<\infty.$$
It follows from Proposition \ref{proposition3.1} that $I_{g_{n-1}}^{n,n-1}:\mathcal{B}^{\alpha}\to\mathcal{B}^{\beta}$ is bounded. Then so is $\sum_{k=0}^{n-2}T_{g_k}^{n,k}$. Through a similar argument, $T_{g_{n-2}}^{n,n-2}:\mathcal{B}^{\alpha}\to\mathcal{B}^{\beta}$ is also bounded. Repeat the above process, we could obtain that each $T_{g_k}^{n,k}:\mathcal{B}^{\alpha}\to\mathcal{B}^{\beta}$ is bounded.    

The proof for the compactness part is similar and we omit the routine details.
\end{proof}

The boundedness and compactness of $I_{g,a}$ and $I_g^{(n)}$ follow immediately as a corollary of Theorem A.

\begin{corollary}
Let $\beta>0$, $\alpha>1$ and $\beta-\alpha+1>0$. Suppose $g\in H(\mathbb{D})$ and $a=(a_0,a_1,\cdots,a_{n-1})\in\mathbb{C}^n$, then $I_{g,a}:\mathcal{B}^{\alpha}\to\mathcal{B}^{\beta}$ is bounded if and only $g\in\mathcal{B}^{\beta-\alpha+1}$. And $I_{g,a}:\mathcal{B}^{\alpha}\to\mathcal{B}^{\beta}$ is compact if and only if $g\in\mathcal{B}^{\beta-\alpha+1}_{0}$.
\end{corollary}

\begin{proof}
Let $g_k=g^{(n-k)}$ in Theorem A, we know that $I_{g,a}:\mathcal{B}^{\alpha}\to\mathcal{B}^{\beta}$ is bounded if and only if $I_{g^{(n-k)}}^{n,k}:\mathcal{B}^{\alpha}\to\mathcal{B}^{\beta}$ is bounded for every $k=0,1,\cdots,n-1$, which is then equivalent to that
\begin{equation}\label{equa3.8}
\sup_{z\in\mathbb{D}}(1-|z|^2)^{\beta-\alpha+n-k}|g^{(n-k)}(z)|<\infty
\end{equation}
for every $k=0,1,\cdots,n-1$ by Proposition 3.1. Since $\beta-\alpha+1>0$, \eqref{equa3.8} is equivalent to that $g\in\mathcal{B}^{\beta-\alpha+1}$ by Lemma 2.1. Similarly, the compactness of $I_{g,a}:\mathcal{B}^{\alpha}\to\mathcal{B}^{\beta}$ is equivalent to that $g\in\mathcal{B}_0^{\beta-\alpha+1}$.
\end{proof}

For $g\in H(\mathbb{D})$, let $I_{g}^{(n)}f=I^n((fg')^{(n-1)})$. When $\alpha>1$, then by Corollary 3.3, the boundedness and compactness of $I_{g}^{(n)}$ are independent of $n$. Specifically, if $I_{g}^{(n)}:\mathcal{B}^{\alpha}\to\mathcal{B}^{\beta}$ is bounded (or compact, resp.) for some $n\geq 1$, then $I_{g}^{(n)}:\mathcal{B}^{\alpha}\to\mathcal{B}^{\beta}$ is bounded (or compact, resp.) for all $n\geq 1$.

We proceed to the proof of Theorem B.

\begin{proof}[\bf Proof of Theorem B]
We only need to prove the necessity. Without loss of generality, assume $N_1\leq N_2\leq \cdots\leq N_n$. If $\sum_{i=1}^{n}I_{g_i}^{N_i,k_i}:\mathcal{B}^{\alpha}\to\mathcal{B}^{\beta}$ is bounded, then
$$\sup_{w\in\mathbb{D}}\left\|\sum_{i=1}^nI_{g_i}^{N_i,k_i}q_w^{[j]}\right\|_{\mathcal{B}^{\beta}}<\infty,\quad j=1,2,\cdots.$$
It follows from Lemma \ref{lemma2.1} that
\begin{equation}\label{equa3.9}
\begin{split}
&\quad \left\|\sum_{i=1}^nI_{g_i}^{N_i,k_i}q_w^{[j]}\right\|_{\mathcal{B}^{\beta}}\\
&\gtrsim (1-|w|^2)^{\beta+N_n-1}\left|\left(\sum_{i=1}^nI_{g_i}^{N_i,k_i}q_w^{[j]}\right)^{(N_n)}(w)\right|\\
&=(1-|w|^2)^{\beta+N_n-1}\left|\sum_{i=0}^n[(q_w^{[j]})^{(k_i)}g_i]^{(N_n-N_i)}(w)\right|.
\end{split}
\end{equation}
Suppose $N_n-N_l+k_l=\max_{1\leq j\leq n}\{N_n-N_j+k_j\}$, then 
$$\max_{\substack{1\leq j\leq n\\j\neq l}}\{N_n-N_j+k_j\}<N_n-N_l+k_l$$
since $N_i-k_i\neq N_j-k_j$ whenever $i\neq j$. Hence by \eqref{equa3.9}, we have
\begin{equation*}
\begin{split}
&\quad\left\|\sum_{i=1}^nI_{g_i}^{N_i,k_i}q_w^{[j]}\right\|_{\mathcal{B}^{\beta}}\\
&\gtrsim (1-|w|^2)^{\beta+N_n+j}\left|(f_w^{[j]})^{(N_n-N_l+k_l)}(w)g_l(w)+\sum_{i=0}^{N_n-N_l+k_l-1}(f_w^{[j]})^{(i)}(w)R_i(w)\right|
\end{split}
\end{equation*}
for $j=1,2,\cdots.$ We modify the proof of Lemma 2.5 and Theorem A to obtain
\begin{equation*}
\sup_{w\in\mathbb{D}}(1-|w|^2)^{\beta-\alpha+N_l-k_l}|g_l(w)|<\infty.
\end{equation*}
Since $k_l>0$, we know that $I_{g_l}^{N_l,k_l}:\mathcal{B}^{\alpha}\to \mathcal{B}^{\beta}$ is bounded by Proposition \ref{proposition3.1}. 

Applying the above process to the boundedness of the sum $\sum_{i\neq l}T_{g_i}^{N_i,k_i}$, we could obtain the boundedness of each $T_{g_i}^{N_i,k_i}$ by induction. The proof for the compactness part is similar, we omit the routine details.
\end{proof}

\begin{remark}
If $\alpha>1$, then Theorem B still holds without the assumption that $k_i>0$ for $i=1,2,\cdots,n$.
\end{remark}

\section{The operator $I_{g}^{n,0}$ in the case $0<\alpha\leq 1$}

\subsection{The case $\alpha=1$}

\begin{proposition}\label{prop4.1}
Let $\beta>0$. Then 
\begin{itemize}
\item[(i)] $I_{g}^{n,0}: \mathcal{B}\to \mathcal{B}^{\beta}$ is bounded if and only 
\begin{equation*}
\sup_{z\in\mathbb{D}}(1-|z|^2)^{n+\beta-1}\log\frac{2}{1-|z|^2}|g(z)|<\infty;
\end{equation*}
\item[(ii)] $I_{g}^{n,0}:\mathcal{B}\to\mathcal{B}^{\beta}$ is compact if and only if 
\begin{equation*}
\lim_{|z|\to 1^-}(1-|z|^2)^{n+\beta-1}\log\frac{2}{1-|z|^2}|g(z)|=0.
\end{equation*}
\end{itemize}
\end{proposition}

\begin{proof}
(a) {\bf Necessity}. For any $w\in\mathbb{D}$, let $h_w(z)=\log\frac{2}{1-\overline{w}z}$, $z\in\mathbb{D}$. It is easy to verify that $\|h_w\|_{\mathcal{B}}\lesssim 1$ for all $w\in\mathbb{D}$. Assume $I_{g}^{n,0}:\mathcal{B}\to\mathcal{B}^{\beta}$ is bounded, then by Lemma \ref{lemma2.1}, we have
\begin{equation*}
\begin{split}
1&\gtrsim\|I_{g}^{n,0}h_w\|_{\mathcal{B}^{\beta}}\\
&\simeq \sup_{z\in\mathbb{D}}(1-|z|^2)^{n+\beta-1}\left| (I_{g}^{n,0}h_w)^{(n)}(z)\right|+\sum_{j=0}^{n-1}|(I_{g}^{n,0}h_w)^{(j)}(0)|\\
&=\sup_{z\in\mathbb{D}}(1-|z|^2)^{n+\beta-1}|h_w(z)||g(z)|\\
&\geq (1-|w|^2)^{n+\beta-1}\log\frac{2}{1-|w|^2}|g(w)|.
\end{split}
\end{equation*}
It follows that
$$\sup_{w\in\mathbb{D}}(1-|w|^2)^{n+\beta-1}\log\frac{2}{(1-|w|^2)}|g(w)|<\infty.$$

{\bf Sufficiency}. For any $f\in \mathcal{B}$, by Lemma \ref{lemma2.2}, $|f(z)|\lesssim \log\frac{2}{1-|z|^2}\|f\|_{\mathcal{B}}$ for all $z\in\mathbb{D}$. Thus, by lemma 2.1,
\begin{equation*}
\begin{split}
\|I_{g}^{n,0}f\|_{\mathcal{B}^{\beta}}&\simeq \sup_{z\in\mathbb{D}}(1-|z|^2)^{n+\beta-1}|I_{g}^{n,0}f(z)|\\
&=\sup_{z\in\mathbb{D}}(1-|z|^2)^{n+\beta-1}|f(z)||g(z)|\\
&\lesssim \|f\|_{\mathcal{B}}\sup_{z\in\mathbb{D}}(1-|z|^2)^{n+\beta-1}\log\frac{2}{1-|z|^2}|g(z)|\lesssim \|f\|_{\mathcal{B}},
\end{split}
\end{equation*}
which implies the boundedness of $I_{g}^{n,0}:\mathcal{B}\to\mathcal{B}^{\beta}$.

(b) {\bf Necessity}. For any $w\in\mathbb{D}$, let 
$$u_w(z)=\frac{(\log\frac{2}{1-\overline{w}z})^2}{\log\frac{2}{1-|w|^2}},\quad z\in\mathbb{D}.$$
It is easy to see that $\sup_{w\in\mathbb{D}}\|u_w\|_{\mathcal{B}}<\infty$ and $\{u_w\}$ converges to 0 uniformly on compact subsets of $\mathbb{D}$ as $|w|\to 1^-$. If $I_{g}^{n,0}:\mathcal{B}\to \mathcal{B}^{\beta}$ is compact, then it follows from Lemma \ref{lemma2.3} that
\begin{equation}\label{equa4.1}
\lim_{|w|\to 1^-}\|I_{g}^{n,0}u_w\|_{\mathcal{B}^{\beta}}=0.
\end{equation}
On the other hand, by Lemma 2.1, we have
\begin{equation*}
\begin{split}
\|I_{g}^{n,0}u_w\|_{\mathcal{B}^{\beta}}&\simeq \sup_{z\in\mathbb{D}}(1-|z|^2)^{n+\beta-1}\left|(I_{g}^{n,0}u_w)^{(n)}(z)\right|\\
&=\sup_{z\in\mathbb{D}}(1-|z|^2)^{n+\beta-1}|u_w(z)||g(z)|\\
&\geq (1-|w|^2)^{n+\beta-1}\log\frac{2}{1-|w|^2}|g(w)|.
\end{split}
\end{equation*}
This, together with \eqref{equa4.1} shows that
$$\lim_{|w|\to 1^-}(1-|w|^2)^{n+\beta-1}\log\frac{2}{1-|w|^2}|g(w)|=0.$$

{\bf Sufficiency}. Let $\{f_j\}$ be any sequence in $\mathcal{B}$ such that $\sup_{j}\|f_j\|_{\mathcal{B}}:=L<\infty$ and $f_j\to 0$ uniformly on compact subsets of $\mathbb{D}$. For any $\varepsilon>0$, choose $0<\delta<1$ such that 
\begin{equation*}
|g(z)|(1-|z|^2)^{\beta+n-1}\log\frac{2}{1-|z|^2}<\frac{\varepsilon}{L}
\end{equation*}
whenever $|z|>\delta$. By Lemma \ref{lemma2.1} and Lemma \ref{lemma2.2}, we have
\begin{equation*}
\begin{split}
\|I_{g}^{n,0}f_j\|_{\mathcal{B}^{\beta}}&\simeq \sup_{z\in\mathbb{D}}(1-|z|^2)^{\beta+n-1}|(I_{g}^{n,0}f_j)^{(n)}(z)|+\sum_{i=0}^{n-1}\left|(I_{g}^{n,0}f_j)^{(i)}(0)\right|\\
&=\sup_{z\in\mathbb{D}}(1-|z|^2)^{\beta+n-1}|f_j(z)||g(z)|\\
&\lesssim \sup_{|z|\leq \delta}|f_j(z)|(1-|z|^2)^{\beta+n-1}|g(z)|\\
&\quad\quad +\sup_{\delta<|z|<1}(1-|z|^2)^{\beta+n-1}|g(z)|\|f_j\|_{\mathcal{B}}\log\frac{2}{1-|z|^2}\\
&\lesssim \sup_{|z|\leq \delta}|f_j(z)|+\varepsilon.
\end{split}
\end{equation*}
Hence $\limsup\limits_{j\to\infty}\|I_{g}^{n,0}f_j\|_{\mathcal{B}^{\beta}}\lesssim\varepsilon$ for any $\varepsilon>0$. It follows from Lemma \ref{lemma2.3} that $I_{g}^{n,0}:\mathcal{B}\to\mathcal{B}^{\beta}$ is compact.
\end{proof}

\subsection{The case $0<\alpha<1$}

\begin{proposition}\label{proposition4.2}
Let $\beta>0$ and $0<\alpha<1$. Then the following conditions are equivalent.
\begin{itemize}
\item[(i)] $I_{g}^{n,0}:\mathcal{B}^{\alpha}\to\mathcal{B}^{\beta}$ is bounded;
\item[(ii)] $I_{g}^{n,0}:\mathcal{B}^{\alpha}\to\mathcal{B}^{\beta}$ is compact;
\item[(iii)] $\sup\limits_{z\in\mathbb{D}}(1-|z|^2)^{n+\beta-1}|g(z)|<\infty$.
\end{itemize}
\end{proposition}

\begin{proof}
It is trivial that (ii) implies (i). By Lemma \ref{lemma2.1}, we have 
$$\sup_{z\in\mathbb{D}}(1-|z|^2)^{n+\beta-1}|g(z)|=\|I_{g}^{n,0}1\|_{\mathcal{B}^{\beta}}.$$
Then (i) implies (ii).

It remains to prove that (iii) implies (ii). So we assume 
$$G:=\sup\limits_{z\in\mathbb{D}}(1-|z|^2)^{n+\beta-1}|g(z)|<\infty.$$
By Lemma \ref{lemma2.1} and Lemma \ref{lemma2.2}, we have
\begin{equation*}
\begin{split}
\|I_{g}^{n,0}f\|_{\mathcal{B}^{\beta}}&\simeq \sup_{z\in\mathbb{D}}(1-|z|^2)^{\beta+n-1}\left|(I_{g}^{n,0}f)^{(n)}(z)\right|+\sum_{j=0}^{n-1}\left|(I_{g}^{n,0}f)^{(j)}(0)\right|\\
&=\sup_{z\in\mathbb{D}}(1-|z|^2)^{\beta+n-1}|g(z)||f(z)|\\
&\lesssim G\|f\|_{\mathcal{B}^{\alpha}}
\end{split}
\end{equation*}
for any $f\in\mathcal{B}^{\alpha}$. This shows that $I_{g}^{n,0}:\mathcal{B}^{\alpha}\to\mathcal{B}^{\beta}$ is bounded.

By Lemma \ref{lemma2.3}, it suffices to prove that $\|I_{g}^{n,0}f_j\|_{\mathcal{B}^{\beta}}\to 0$ for any bounded sequence $\{f_j\}$ in $\mathcal{B}^{\alpha}$ which converges to 0 uniformly on compact subsets of $\mathbb{D}$. To this end, let $L:=\sup_{j}\|f_j\|_{\mathcal{B}^{\alpha}}$. For any $\varepsilon>0$, choose $\delta\in (\frac{1}{2},1)$ such that $(1-\delta)^{1-\alpha}<\frac{1-\alpha}{L}\varepsilon$. Then 
\begin{equation}\label{equa4.2}
\begin{split}
\left|f_j(z)-f_j\left(\frac{\delta}{|z|}z\right)\right|&=\left|\int_{\frac{\delta}{|z|}}^1zf_j'(tz)dt\right|\leq \int_{\frac{\delta}{|z|}}^1\frac{|z|}{(1-|tz|^2)^{\alpha}}dt\|f_j\|_{\mathcal{B}^{\alpha}}\\
&\lesssim L\int_{\delta}^1\frac{1}{(1-x)^{\alpha}}dx=\frac{L(1-\delta)^{1-\alpha}}{1-\alpha}<\varepsilon.
\end{split}
\end{equation}
Using Lemma \ref{lemma2.1} again, we have
\begin{equation*}
\begin{split}
\|I_{g}^{n,0}f_j\|_{\mathcal{B}^{\beta}}&\simeq \sup_{z\in\mathbb{D}}(1-|z|^2)^{\beta+n-1}|g(z)||f_j(z)|\\
&\leq G\left(\sup_{|z|\leq \delta}|f_j(z)|+\sup_{\delta<|z|<1}|f_j(z)|\right).
\end{split}
\end{equation*}
By assumption, we know that 
$$\lim_{j\to\infty}\sup_{|z|\leq \delta}|f_j(z)|=0.$$
On the other hand, by triangle inequality and \eqref{equa4.2}, we get
\begin{equation*}
\begin{split}
\sup_{\delta<|z|<1}|f_j(z)|&\leq \sup_{\delta<|z|<1}\left|f_j(z)-f_j\left(\frac{\delta}{|z|}z\right)\right|+\sup_{\delta<|z|<1}\left|f_j\left(\frac{\delta}{|z|}z\right)\right|\\
&<\varepsilon+\sup_{|z|=\delta}|f(z)|
\end{split}
\end{equation*}
Therefore, 
$$\lim_{j\to\infty}\sup_{\delta<|z|<1}|f_j(z)|<\varepsilon+\lim_{j\to\infty}\sup_{|z|=\delta}|f_j(z)|=\varepsilon.$$
It follows that $\limsup\limits_{j\to\infty}\|I_{g}^{n,0}f_j\|_{\mathcal{B}^{\beta}}<G\varepsilon$ for any $\varepsilon>0$. Hence $\|I_{g}^{n,0}f_j\|_{\mathcal{B}^{\beta}}\to 0$, and the proof is complete.
\end{proof}

\subsection{The solution of linear differentiation equations}

As an application, we can investigate the solution of certain linear differentiation equations. For a vector-valued analytic function $\mathbf{g}=(g_0,g_1,\cdots,g_{n-1})$, let 
$$\|\mathbf{g}\|_*=\max\{\|g_k\|_{k}:0\leq k\leq n-1\},$$
here, for $g\in H(\mathbb{D})$, $\|g\|_{k}=\sup_{z\in\mathbb{D}}(1-|z|^2)^{n-k}|g(z)|$ when $1\leq k\leq n-1$ and 
\begin{equation*}
\|g\|_0=\left\{
  \begin{array}{lr}
\sup\limits_{z\in\mathbb{D}}(1-|z|^2)^{n}|g(z)|,&\quad \alpha>1,\\
\sup\limits_{z\in\mathbb{D}}(1-|z|^2)^{n}\log\frac{2}{1-|z|^2}|g(z)|, &\quad \alpha=1,\\
\sup\limits_{z\in\mathbb{D}}(1-|z|^2)^{n+\alpha-1}|g(z)|, & \quad 0<\alpha<1.
  \end{array}
\right. 
\end{equation*}

\begin{proposition}
Let $\alpha>0$, $F\in\mathcal{B}^{\alpha+n}$ and $\mathbf{g}=(g_0,g_1,\cdots,g_{n-1})$ be an analytic vector-valued function. Then there exists an $A>0$ such that every solution of the following nonhomogeneous linear differential equation
$$f^{(n)}+f^{(n-1)}g_{n-1}+\cdots+fg_0=F$$
belongs to $\mathcal{B}^{\alpha}$ whenever $\|\mathbf{g}\|_{*}<A$.
\end{proposition}

\begin{proof}
If $f$ is any solution of the above differential equation, let 
$$F_0(z)=(I^nF)(z)+\sum_{j=0}^{n-1}\frac{1}{j!}f^{(j)}(0)z^j,\quad z\in\mathbb{D}.$$
By Lemma \ref{lemma2.1}, we know that $I^nF\in\mathcal{B}^{\alpha}$ since $F\in\mathcal{B}^{\alpha+n}$. It follows that $F_0\in \mathcal{B}^{\alpha}$ and $F_0^{(j)}(0)=f^{(j)}(0)$ for $0\leq j\leq n-1$. 
Recall the operator $I_{\mathbf{g}}^{(n)}$ defined by
$$I_{\mathbf{g}}^{(n)}f=I^n(fg_0+f'g_1+\cdots+f^{(n-1)}g_{n-1}).$$
By Propositions 3.1, 4.1 and 4.2, we know that $I_{\mathbf{g}}^{(n)}$ is bounded on $\mathcal{B}^{\alpha}$, and there exists a constant $C>0$, depending only on $\alpha$, such that
\begin{equation*}
\|I_{\mathbf{g}}^{(n)}\|\leq C\sum_{k=0}^{n-1}\|g_k\|_{k}.
\end{equation*}
Choose $A<\frac{1}{2nC}$, then the norm $\|I_{\mathbf{g}}^{(n)}\|<1$, and hence there exists $\tilde{f}\in\mathcal{B}^{\alpha}$ such that $I_{\mathbf{g}}^{(n)}\tilde{f}+\tilde{f}=F_0$. It is easy to see that $\tilde{f}$ and $f$ satisfy the same differential equation with the same initial condition. By the uniqueness of the solution of the initial value problem, $f=\tilde{f}\in\mathcal{B}^{\alpha}$.
\end{proof}

\begin{remark}
We note that if $(1-|z|^2)^{n-k}|g_k(z)|\to 0$ as $|z|\to 1^-$ for $0\leq k\leq n-1$, and 
\begin{equation*}
\left\{
  \begin{array}{lr}
(1-|z|^2)^{n}|g(z)|\to 0,&\quad \alpha>1,\\
(1-|z|^2)^{n}\log\frac{2}{1-|z|^2}|g(z)|\to 0, &\quad \alpha=1,\\
(1-|z|^2)^{n+\alpha-1}|g(z)|\to 0, & \quad 0<\alpha<1,
  \end{array}
\right. 
\end{equation*}
then the same result in Proposition 4.3 holds without any restriction on the norm $\|\mathbf{g}\|_{*}$. This comes from the fact that $I_{\mathbf{g}}^{(n)}$ is compact on $\mathcal{B}^{\alpha}$ if the above conditions are satified. In this case, we claim that the spectrum of $I_{\mathbf{g}}^{(n)}$ is the singleton $\{0\}$. Then the desired result will follow from the same arguments above. In fact, if $\lambda\in \mathbb{C}\backslash\{0\}$, then the equation $I_{\mathbf{g}}^{(n)}f=\lambda f$ is equivalent to the initial value problem
$$f^{(n-1)}g_{n-1}+\cdots+f'g_1+fg_0=\lambda f^{(n)},\quad f(0)=\cdots=f^{(n-1)}(0)=0.$$
By the uniqueness of solutions, we must have $f=0$. Thus $I_{\mathbf{g}}^{(n)}$ has no nonzero eigenvalue.
\end{remark}

\begin{declaration}
 The authors declare that there are no conflicts of interest regarding the publication of this paper. No data was generated by this project.
\end{declaration}

\end{document}